\documentclass[10pt]{article}  
\usepackage{amsmath}
\usepackage{amssymb}
\usepackage{epsfig}
\usepackage{changebar}
\usepackage{theorem}

\theoremheaderfont{\scshape}
\newtheorem{thm}{Theorem}[section]
\newtheorem{cor}[thm]{Corollary}
\newtheorem{lem}[thm]{Lemma}

\theorembodyfont{\normalfont}
\newtheorem{rem}{Remark}[section]

\newenvironment{proof}{

\noindent{\it Proof. }}
{\hfill\rule{2mm}{2mm}\vskip3mm \par} 
\numberwithin{equation}{section}

\newcommand{\ds}{\displaystyle}

\newcommand{\fracp}[2]{\frac{\partial #1}{\partial #2}}
\newcommand{\fracpsec}[2]{\frac{{\partial^2 #1}} {{\partial #2^2}}}
\newcommand{\fracpsecc}[3]{\frac{{\partial^2 #1}}{{\partial #2 \partial #3}}}
\newcommand{\fracd}[2]{\frac{d #1}{d #2}}

\newcommand{\nit}{\mathbb{N}}
\newcommand{\rit}{\mathbb{R}}

\newcommand{\HH}{\mathbf{H}}
\newcommand{\EE}{\mathbf{E}}
\title{Existence result for a model of {\it Proteus mirabilis} 
swarm\footnote{This work has been partially supported by ``La R\'egion Bretagne'',
F-35031 Rennes, Program 1042: ``Renouvellement des comp\'etences dans les laboratoires 
de recherche'', Operation A1C872. }
}
\author{
E. Fr\'enod\thanks{LMAM et Lemel, 
Universit\'e de Bretagne Sud, Centre Yves Coppens, Campus de Tohannic,
F-56000, Vannes.}
}

\begin{document}
\maketitle

{\small {\bf Abstract: } In this paper we present a modification of the usual 
{\it Proteus mirabilis} Swarm model. For the obtained model
(which is a two phase model with a non-linear diffusion term containing memory)
we set up a collection of a priori estimates. 
Those estimates allow to get an existence and uniqueness result.}
\section{Introduction and results}\label{SEcIntRes}
{\it Proteus mirabilis} is a bacterium that can be either a short cell
we call ``swimmer'' or an elongated cell capable of translocation
we call ``swarmer''.
A model of behaviour of {\it Proteus mirabilis} colony has been proposed
by Esipov and Shapiro \cite{EsiSha}, based on ideas of Gurtin \cite{Gur}.
\newline
In this paper, we prove an existence result to a model which is,
in a way, a generalization but also a regularization of the
Esipov and Shapiro \cite{EsiSha} model. 

~

The model under consideration here is a two phase model with
a non-linear diffusion term containing memory for one of 
the two phases.
It involves two functions $\rho$ and
$Q$. The function $\rho=\rho(t,a,x)$ refers, at time $t\in[0,T)$,
$0<T<+\infty$, to the density of swarmers of age $a\in [0,A)$, 
$0<A\leq +\infty$ at position $x\in\Omega$, where $\Omega$ is a regular
sub-domain of $\rit^2$, with boundary $\partial\Omega$. In each point $x$ of
$\partial\Omega$, $\overrightarrow{\nu}= \overrightarrow{\nu}(x)$ stands for the
unit  normal vector pointing outside $\Omega$. The function $Q=Q(t,x)$
stands, at time $t$, for the biomass density of swimmers on $\Omega$.
\newline
For a constant $\tau$, those
two functions are supposed to satisfy the following system:
\begin{align}
& \frac{\partial Q}{\partial t}=\frac{1-\xi}{\tau}Q+
\int_0^{A}\rho(\cdot ,a,\cdot)e^{a/\tau}\mu(\cdot,a,\cdot)\,da +
\chi(A)\rho(\cdot,A,\cdot)e^{A/\tau},
~~\text{ on } [0,T) \times \Omega,
\label{M1.a}
\\
&\frac{\partial\rho}{\partial t}+\frac{\partial\rho}{\partial a}=
-\mu\rho+
\nabla\cdot\big[\big(D({\cal M},Q,P)+d\big)\nabla\rho\big],
~~\text{ on } [0,T) \times[0,A)\times \Omega,
\label{M1.b}
\\
&\rho(\cdot,0,\cdot)= \frac{\xi}{\tau}Q,
~~\text{ on } [0,T) \times \Omega,
\label{M1.c}
\\
&\rho(0,\cdot,\cdot)= \rho_0,
~~\text{ on } [0,A) \times \Omega,
\label{M1.d}
\\
& Q(0,\cdot)=Q_0,
~~\text{ on } \Omega,
\label{M1.e}
\\
&\fracp{\rho}{\overrightarrow{\nu}} = 0,
~~\text{ on } [0,T) \times[0,A) \times \partial\Omega.
\label{M1.f}
\end{align}
Above, $\chi(A)$ is an artifice allowing to take into account a possible
maximum age $A$ beyond which swarmers cannot exist. It has the 
following definition:
\begin{equation}\label{Ki}
\chi(A)=1,~~ \forall A\in \rit,~~~~~~~~~\chi(+\infty)=0.
\end{equation}
Denoting by $C^k_b$ the space of functions having continuous and bounded 
derivatives up to order $k$, $\mu=\mu(t,a,x)$ is a function such that
\begin{align}\label{Hypmu}
& \mu\in C^2_b([0,T) \times[0,A)\times \Omega),~~ \mu\geq 0,
~~\lim_{a\rightarrow A}\mu(t,a,x) =\overline{\mu}
\text{ ~ uniformly in $x$ and $t$,}
\end{align}  
with $\overline{\mu}\geq c(1-\chi(A))$, for a constant $c>0$.
The function $\xi=\xi(t,Q)$ satisfies
\begin{equation}\label{Hypxi}
\xi\in C^2_b([0,T) \times \rit),~~~~ 0\leq \xi\leq 1.
\end{equation}
In the second equation, $\nabla$ stands for the gradient with respect to the
$x-$variable, and  $\nabla\cdot$ for the divergence. The diffusion coefficient
is the sum of a constant, a priori small,
\begin{equation}\label{HypD}
d>0,
\text{ ~ and of ~ } 
D=D({\cal M},P,Q), \text{ a non negative } C^1_b \text{ function of its arguments.}
\end{equation} 
In (\ref{HypD}), $P$ is defined for $0\leq a_{min} < A$ by
\begin{equation}\label{DefP}
P(t,x) = \int_{a_{min}}^A \rho(t,a,x) e^{a/\tau} da,
\end{equation}
and $Q$ is given by the first equation of the system.
\newline
The memory (or hysteresis) term ${\cal M} = {\cal M}_{[P]}(t,x)$ keeps
information on the value of $P$ in the past. For four thresholds,
$P_{min}<p_{min}<p_{max}<P_{max}$, with $P_{min}$ close to $p_{min}$ and $P_{max}$
close to $p_{max}$, ${\cal M}$ is defined as the solution to:
\begin{align}
&\fracp{\cal M}{t} = \frac{1}{P_{max}-p_{max}} 
  ~H_r\bigg(\frac{P-p_{max}}{P_{max}-p_{max}}\bigg)
  ~H_r (1- {\cal M})-  \frac{1}{p_{min}-P_{min}}
  ~H_r\bigg(\frac{p_{min}-P}{p_{min}-P_{min}}\bigg)~H_r ({\cal M}),
\nonumber
\\
& {\cal M}(0,\cdot)= {\cal M}_0,
\label{Defm}
\end{align}
with, denoting $P_0=P(0,\cdot)$,
\begin{equation}\label{Hypm0}
{\cal M}_0 \in C^1_b(\Omega), ~~ 0\leq {\cal M}_0 \leq 1, ~~ 
{\cal M}_0=0 \text{ where } P_0<P_{min} \text{ and }
{\cal M}_0=1 \text{ where } P_0>P_{max},
\end{equation}\label{DefrH}
and with
\begin{equation}
H_r(p) =0 \text{ if } p\leq 0,  ~~ H_r(p) =p \text{ if } 0\leq p \leq 1 
 ~~\text{ and } H_r(p) =1 \text{ if } p\geq 1.
\end{equation}

~

We now turn to the statement of the main result of this paper.
\begin{thm}\label{Thm1}
Under assumptions (\ref{Hypmu}), (\ref{Hypxi}), (\ref{HypD}) and 
(\ref{Hypm0}), if 
$\rho_0\geq 0 \in L^1 \cap W^{2,2} \cap W^{1,4}
([0,A)\times\Omega,~e^{a/\tau}da dx)$ satisfies
\begin{equation}\label{Hyprho0}
\begin{aligned}
&\int_\Omega  (\rho_0)^2 e^{2a/\tau} dx \leq
b \int_\Omega  (\rho_0)^2 e^{a/\tau} dx,  ~~\forall a\in [0,A),
\\
&\int_0^A \int_\Omega |\nabla\rho_0|^4 e^{4a/\tau}da dx \leq
b \int_0^A \int_\Omega |\nabla\rho_0|^4 e^{a/\tau}da dx, 
\\
&\int_0^A \int_\Omega |\nabla\rho_0|^2 e^{2a/\tau} dx \leq
b \int_0^A\int_\Omega |\nabla\rho_0|^2 e^{a/\tau} dx,
\end{aligned}
\end{equation}
for a constant $b$, and if $Q_0 \geq 0 \in L^1 \cap W^{2,2} \cap W^{1,4}(\Omega)$;
then, there exists a unique solution 
$(Q,\rho)\in L^\infty\big(0,T; ((L^1 \cap W^{1,2}(\Omega) )
\times (L^1 \cap W^{1,2}([0,A)\times\Omega,~e^{a/\tau}da dx)) \big)$
to system (\ref{M1.a})-(\ref{M1.f}) coupled with (\ref{DefP}) and 
(\ref{Defm}). Moreover, $Q\geq 0$ and $\rho\geq 0$.
\end{thm}
The precise definitions of the spaces at work in the Theorem are given in
the beginning of section \ref{SecExUn}.

~

We now give references where modelling and mathematical methods are
developed on age-structured population problem:
Gurtin and Mac Camy \cite{GurMcCa1974},
Marcati \cite{Mar1981},
Andreasen \cite{Andr1989,Andr1992,Andr1995}; 
possibly with diffusion:
Gurtin \cite{Gur},
Di Blasio and Lamberti \cite{DiBlaLam1978},
Di Blasio \cite{DiBla1979},
Mac Camy  \cite{McCa},
Gurtin and Mac Camy \cite{GurMcCa}, 
Busenberg and  Iannelli \cite{BusIann1983},
Langlais \cite{Langlais1985,Langlais1988},
Kubo and Langlais \cite{KuLang1991},
Huang \cite{Hua1994}
and 
Esipov and Shapiro \cite{EsiSha}.
For simulation methods we refer to
Lopez and Trigiante \cite{LopTri1985},
Milner \cite{Mil1990}
Kim \cite{Kim1996},
Esipov and Shapiro \cite{EsiSha},
Medvedev, Kapper and Koppel \cite{MedKaKo},
Ayati and Dupont \cite{AyDu}.
Concerning the biological description of {\it Proteus mirabilis}
colony behaviour, we refer for instance to 
Rauprich {\em et al} \cite{Rauetal1996},
Gu\'e, Dupont, Dufour and Sire \cite{Gueetal2001}
and theirs references. 

~

The paper is organized as follows: in section \ref{SecMod} we present
the way to go from the Esipov and Shapiro model to system 
(\ref{M1.a})-(\ref{M1.f}). Then section \ref{SecAPrEst} is devoted 
to a priori estimates for the solution to (\ref{M1.a})-(\ref{M1.f}).
By a usual procedure consisting in linearizing and passing to limit, we
prove the Theorem in  section \ref{SecExUn}.

~
\newline
{\bf Acknowledgements:} I would like to thank O. Sire for having introduced
to me the swarm model of Esipov and Shapiro \cite{EsiSha} and for stimulating
discussions.
\newline
I would also like to thank F. Granger who, despite stopping its PhD thesis
for personal reasons, made the first steps towards the result.

\section{Model}\label{SecMod}
{\it Proteus mirabilis} is a pathogenic bacterium of urinary tract that
when standing in liquid medium, consists in a usual short 
``swimmer cell'' or ``swimmer''. When placed on agar medium, if
the bacterium density is large enough, it may undergo a differentiation
process producing an elongated cell with several nuclei called 
``swarmer cell'' or ``swarmer''. Those swarmers are capable of translocation
allowing the bacterial colony to colonise the medium.

 ~

The macroscopic model built by Esipov and Shapiro \cite{EsiSha}
describes this swarm phenomenon at the colony scale. We shall explain
this model now. The swarmer behaviour depends on their own age.
This dependence is taken into account by introducing the age dependent
density of swarmers $\rho(t,a,x)$. The link between this density $\rho$
and the biomass density is a consequence of the fact that 
the mass of each cell is in direct proportion with its length
and that the length growth of a swarmer is also in direct proportion
with its length. Then the biomass density, at time $t$, of swarmers
of age $a$ at position $x$ is $\rho(t,a,x)e^{a/\tau}$, where 
$\tau$ is the growth rate of the biomass. The first age-depending 
behaviour of the swarmers is that they actively participate in group
migration only after an age $a_{min}$. Then the definition of the
biomass density $P$ of swarmers capable of active translocation
is given by (\ref{DefP}). The second age-depending behaviour is that
the swarmers dedifferentiate themselves and give swimmers. On this topic, 
Esipov and Shapiro \cite{EsiSha} consider two situations. In the
first one (Model A), the swarmers dedifferentiate themselves at a given 
age $a_{max}$. The second situation consider that swarmers may 
dedifferentiate at each time with a probability $1/\overline{a}$
(Model B).
\newline
Now, we are able to write the swimmer evolution equation for the 
biomass density $Q$. Its evolution results from the classical cellular 
division, with a characteristic time which is the biomass growth rate 
$\tau$, subtracting the proportion of bacteria undergoing 
differentiation and adding the dedifferentiation product. 
In the case of Model A, the evolution equation for $Q$ is then
\begin{equation}\label{EqQA}
\fracp{Q}{t}=\frac{1-\xi}{\tau} Q + 
\rho(\cdot,a_{max},\cdot) e^{a_{max}/\tau},
~~~~ Q(0,\cdot)=Q_0.
\end{equation}
Here, $Q_0$ stands for the initial swimmer density and $\xi/\tau$
is the fraction of swimmer population to produce swarmers.
In the case of Model B, the evolution equation is
\begin{equation}\label{EqQB}
\fracp{Q}{t}=\frac{1-\xi}{\tau} Q + 
\int_0^t
\rho(\cdot,a,\cdot) \,\frac{e^{a/\tau}}{\overline{a}} da,
~~~~ Q(0,\cdot)=Q_0.
\end{equation}
\newline
We turn to the evolution of the swarmer density $\rho$. Its evolution is
linked with ageing and the dedifferentiation process, but also to swarm.
This last phenomenon is modelled by a non linear diffusion term with 
memory. The evolution equation for $\rho$ is then
\begin{equation}\label{EqrhoA}
\frac{\partial\rho}{\partial t}+\frac{\partial\rho}{\partial a}=
\nabla\cdot\big[D({\cal M},Q,P)\nabla\rho\big],
\end{equation}
in the case of Model A; and, in the case of Model B, it is  
\begin{equation}\label{EqrhoB}
\frac{\partial\rho}{\partial t}+\frac{\partial\rho}{\partial a}=
-\frac{1}{\overline{a}}\rho+
\nabla\cdot\big[D({\cal M},Q,P)\nabla\rho\big].
\end{equation}
Both of those equations are equiped with the following initial
and boundary conditions:
\begin{equation}\label{IBCond}
\text{(a): } \rho(\cdot,0,\cdot)= \frac{\xi}{\tau}Q,~~~~~~~~
\text{(b): } \rho(0,\cdot,\cdot)= 0,~~~~~~ ~~
\text{~~~~ and ~~~~~~~} 
\text{(c): } {\fracp{\rho}{\overrightarrow{\nu}}}_{\big|\partial\Omega} = 0. 
\end{equation}
The first of those three conditions means the fraction $\xi/\tau$ of 
swimmers undergoing the differentiation process produces swarmers of age 0.
The initial condition  on $\rho$ means that, at the beginning of the process,
there is no swarmer. The boundary condition means that swarmer cannot leave
the domain $\Omega$.
In the diffusion term 
$\nabla\cdot\big[D({\cal M},Q,P)\nabla\rho\big]$
appearing in (\ref{EqrhoA}) and (\ref{EqrhoB}), and modelling the
swarm, the diffusion factor $D({\cal M},Q,P)$ depends on the present 
value of $Q$ and $P$ but also on the history of $P$.
\newline
The term ${\cal M}$ then keeps in memory informations concerning
the history of the swarmer density.  Esipov and Shapiro \cite{EsiSha}
defines ${\cal M}$ as being set to 1 in a given point $x$ if the value 
of $P$ in $x$ reaches a threshold $P_{max}$. Then it remains at the
value 1 until the value of $P$ in $x$ reaches another value
$P_{min}<P_{max}$. Then they suggest to take
\begin{equation}\label{DefD}
D({\cal M},Q,P) = \overline{D_0} \, {\cal M} \; \gamma\Big(\frac{P}{P_{max}}\Big)
\;\exp\Big( \frac{-Q}{Q_{sat}} \Big),
\end{equation}
for given values of $\overline{D_0}$ and $Q_{sat}$ and with,
$H$ being such that $H(p)=0$ if $p<0$ and $H(p)=1$ if $p>0$,
\begin{equation}\label{Defgamma}
\gamma(p) = (p-\frac{P_{min}}{P_{max}}) H(p-\frac{P_{min}}{P_{max}})\text{ or }
\gamma(p) = p-\frac{P_{min}}{P_{max}} \text{ or }
\gamma(p) = p^2  \text{ or } \gamma(p) = 1.
\end{equation}
We also mention that Medvedev, Kapper and Koppel \cite{MedKaKo}
studied this model taking, for a given value of $k$,
\begin{equation}\label{DefDMKK}
D({\cal M},Q,P) = D(Q,P) = \frac{\overline{D_0} P}{P+kQ}.
\end{equation}

~

Now, we explain in what sense the model (\ref{M1.a})-(\ref{M1.f}) is a 
generalization and a regularization of the Esipov and Shapiro \cite{EsiSha}
model.
\newline
First, because of the initial condition (\ref{IBCond}.b), it is an easy
game to see that, at least formally, the solution $\rho$ to
(\ref{EqrhoA}) or (\ref{EqrhoB}) satisfies 
\begin{equation}
\rho(t,a,x) =0 \text{ when } a>t.
\end{equation}
Hence we can replace (\ref{EqQB}) by
\begin{equation}\label{EqQB2}
\fracp{Q}{t}=\frac{1-\xi}{\tau} Q + 
\int_0^{+\infty}
\rho(\cdot,a,\cdot) \,\frac{e^{a/\tau}}{\overline{a}} da,
~~~~ Q(0,\cdot)=Q_0.
\end{equation}
Making this allows to take under consideration, with no loss of 
consistency, initial data $\rho(0,\cdot,\cdot)$ that are not 0
coming to the more general initial and boundary conditions
(\ref{M1.c})-(\ref{M1.f}).
\newline
Secondly, equations (\ref{EqQA}) and (\ref{EqQB2}) are particular
cases of the general equation  (\ref{M1.a}) with assumption
(\ref{Hypmu}). The case of equation (\ref{EqQA}) is recovered 
setting $\mu=0$ and $A=a_{max}$ and the case of (\ref{EqQB2})
is recovered setting $A=+\infty$ and $\mu=1/\overline{a}$.
We also see that (\ref{M1.b}) with $d=0$ is a general framework
inside which (\ref{EqrhoA}) and (\ref{EqrhoB}) may enter directly.
\newline
Concerning the fraction $\xi$ of swimmers to produce swarmers,
it seems to depend on experimental conditions and to be 0 when
the swimmer density is high. It seems then to be reasonable to
set that $\xi$ satisfies (\ref{Hypxi}). 
\newline
The first regularization effect we consider in our model 
(\ref{M1.a})-(\ref{M1.f}) consists in adding $d>0$ in 
(\ref{M1.b}). This may be justified by experimental arguments
saying that swarmers always experience a small random motion even 
before and after swarming. 
\newline
The second regularization effect which constitutes the most 
visible modification of the model concerns the memory term
${\cal M}$. In order to explain the way to go from the 
definition of ${\cal M}$ by  Esipov and Shapiro \cite{EsiSha} to
the definition (\ref{Defm}), we first notice that the term ${\cal M}$ 
of Esipov and Shapiro could be defined formally, in any $x$,
 as the solution to the following equation:
\begin{equation}\label{DefMES}
\fracp{{\cal M}}{t}(\cdot,x) = (1-{\cal M}(\cdot,x))
\delta_{\{t / P(t,x)=P_{max}\}}
- {\cal M}(\cdot,x) \delta_{\{t / P(t,x)=P_{min}\}},
~~~ {\cal M}(0,x)=0,
\end{equation}
where $\delta_{\{t / P(t,x)=P_{max}\}}$ stands for the Dirac 
measure in instant $t$ where $P(t,x)=P_{max}$. Of course
this equation has no real mathematical meaning. But formally,
if ${\cal M}=0$ and if the value $P_{max}$ is reached at a given
time $\tilde t$, the solution of (\ref{DefMES}) experiences
a jump $+1$. When $P=P_{max}$, and ${\cal M}=1$, nothing happens
and this is what is needed. In the same way, if at a given
time $\tilde t$,  $P=P_{min}$, nothing happens if ${\cal M}=0$
and ${\cal M}$ experiences a jump $-1$ if ${\cal M}=1$.
Now, it is easy to see that the right hang side of (\ref{Defm}) is
nothing but a regularization of the right hand value of (\ref{DefMES}).
This is the reason why we make this choice to define ${\cal M}$.

\section{A priori estimates}\label{SecAPrEst}
The key point to get the existence result is a collection
of a priori estimates satisfied by $(Q,\rho)$.
They are mathematical translations of biological properties.
In order to set those estimates comfortably, we assume that
the assumptions and the conclusions of Theorem \ref{Thm1}
are satisfied.
\subsection{$L^1$ and $L^2$ estimates}\label{SecL1L2Es}
For $p>0,$ we denote
\begin{equation}\label{DefNorm}
\| \rho(t)\|_{p}=
\bigg(\int_0^A \int_\Omega |\rho(t,a,x)|^p \,e^{a/\tau} dxda\bigg)^{1/p},~~
\| Q(t)\|_{p}=\bigg(\int_\Omega |Q(t,x)|^p \, dx\bigg)^{1/p},
\end{equation}
and we have the following estimates saying that the total biomass
grows exponentially with a growth rate $\tau$.
\begin{lem}\label{LemEs1}
If the assumptions of Theorem \ref{Thm1} are valid, and if 
the solution $(Q,\rho)$ given by this same Theorem exists,
then it satisfies
\begin{equation}\label{EsL1}
\| \rho(t)\|_{1} + \| Q(t)\|_{1} = 
\big(\| \rho_0\|_{1} + \| Q_0\|_{1}\big)e^{t/\tau}.
\end{equation}
\end{lem}

\begin{proof}
Multiplying equation (\ref{M1.b}) by $e^{a/\tau}$, integrating then
with respect to $x$ and $a$ yields:
\begin{equation}
\fracd{\| \rho \|_{1}}{t} +
\lim_{a\rightarrow A} \bigg(\int_\Omega \rho\,e^{a/\tau}dx\bigg)-
\int_\Omega \rho(\cdot,0,\cdot)\,dx = 
\frac1\tau \| \rho \|_{1} -
\int_0^A \int_\Omega \mu \rho \,e^{a/\tau}dxda,
\end{equation}
which also reads 
\begin{equation}\label{Es1.1}
\fracd{\| \rho \|_{1}}{t} +
\chi(A) e^{A/\tau}\int_\Omega \rho(\cdot,A,\cdot)\,dx-
\int_\Omega \frac{\xi}{\tau} Q \,dx = 
\frac1\tau \| \rho \|_{1} -
\int_0^A \int_\Omega \mu \rho \,e^{a/\tau}dxda.
\end{equation}
Integrating now (\ref{M1.a}) with respect to $x$, we obtain
\begin{equation}\label{Es1.2}
\fracd{\| Q \|_{1}}{t} = \int_\Omega \frac{1-\xi}{\tau} Q \, dx+
\int_0^A \int_\Omega \mu \rho \,e^{a/\tau}dxda + 
\chi(A) e^{A/\tau}\int_\Omega \rho(\cdot,A,\cdot)\,dx.
\end{equation}
Summing up (\ref{Es1.1}) and (\ref{Es1.2}) finally gives
\begin{equation}\label{Es1.3}
\fracd{\big(\| \rho \|_{1} +\| Q \|_{1}\big)}{t} = 
\frac1\tau \big(\| \rho\|_{1} + \| Q \|_{1} \big),
\end{equation}
proving the Lemma. ~~~~~
\end{proof}
The second Lemma concerns the $L^2-$norms of $\rho$ and $Q$.
It mathematically translates that the biomass cannot be so
gathered that a null area set contains a positive biomass
quantity.
\begin{lem}\label{LemEs2}
If the assumptions of Theorem \ref{Thm1} are realized and if the 
solution $(Q,\rho)$ exists, then for any $t\in [0,T)$, it satisfies
\begin{equation}\label{EsL2}
\| \rho(t) \|_{2} + \| Q(t) \|_{2} \leq  
c\big( \| \rho_0 \|_{2} + \| Q_0 \|_{2} \big),
\end{equation} 
for a constant $c$ ($c$ only depends on $A$, $\tau$, $T$ and  
$\tilde \mu =\sup \{\mu(t,a,x),\, t\in [0,T), a\in [0,A), x\in\Omega \}.$)
\end{lem}
\begin{proof}
First, integrating (\ref{M1.b}) with respect to $x$, we get:
\begin{equation}\label{EsL2.1}
\bigg(\fracp{}{t}+\fracp{}{a}\bigg) 
\bigg( \int_\Omega \rho \, dx \bigg)=
- \int_\Omega \mu\rho \, dx \leq 0.
\end{equation}
Making the same, after multiplying (\ref{M1.b}) by $\rho$, gives:
\begin{equation}\label{EsL2.2}
\frac12 \bigg(\fracp{}{t}+\fracp{}{a}\bigg)
\bigg( \int_\Omega \rho^2 \, dx \bigg)=
- \int_\Omega \mu\rho^2 \, dx -
\int_\Omega (D+d) |\nabla\rho|^2 \,dx\leq 0.
\end{equation}
Defining, for a fixed $\alpha$ and for $p=1$ or $2$,
\begin{equation}\label{EsL2.3}
r_p:s\mapsto \int_\Omega \rho^p(s,\alpha+s,x) \, dx ,
\end{equation}
we have 
\begin{equation}\label{EsL2.4}
r_p'(s) = \Bigg(\bigg(\fracp{}{t}+\fracp{}{a}\bigg) 
\bigg(\int_\Omega \rho^p\, dx \bigg) \Bigg)(s,\alpha+s) \leq 0.
\end{equation}
Now, since for fixed $t$ and $a$ and with $\alpha=a-t$ we have 
$\int_\Omega \rho^p(t,a,x)\, dx =r_p(t)$, if $t<a$ the relation 
$r_p(t)\leq r_p(0)=  \int_\Omega \rho^p(0,a-t,x)\, dx$ reads
\begin{equation}\label{EsL2.5}
\int_\Omega \rho(t,a,x)\, dx \leq \int_\Omega \rho_0(a-t,x)\, dx , ~~
\int_\Omega \rho^2(t,a,x)\, dx \leq \int_\Omega \rho_0^2(a-t,x)\, dx.
\end{equation}
In the case when $t>a$, the relation 
$r_p(t)\leq r_p(t-a)= \int_\Omega \rho^p(t-a,0,x)\, dx$ gives
\begin{equation}\label{EsL2.6}
\int_\Omega \rho(t,a,x)\, dx \leq 
\int_\Omega\frac{\xi}{\tau} Q(t-a,x)\, dxda,~~
\int_\Omega \rho^2(t,a,x)\, dx \leq 
\int_\Omega\frac{\xi^2}{\tau^2}  Q^2(t-a,x)\, dxda,
\end{equation}
Secondly, we multiply equation (\ref{M1.b}) by $e^{a/\tau}\rho$ and
we integrate in $x$ and $a$. Since
\begin{equation}\label{EsL2.7}
\int_0^A\fracp{\rho}{a}\rho \, e^{a/\tau}da =
-\int _0^A\fracp{\rho}{a}\rho \, e^{a/\tau}da-
\frac{1}{\tau} \int _0^A\rho^2 \, e^{a/\tau}da 
+\bigg[\rho^2 \, e^{a/\tau} \bigg]_{a=0}^{a=A},
\end{equation}
and since 
\begin{equation}\label{EsL2.8}
\int_\Omega \rho^2 (t,0,x)\, dx = 
\int_\Omega \frac{\xi^2}{\tau^2} Q^2 \, dx,
\end{equation}
we obtain
\begin{equation}\label{EsL2.9}\begin{aligned}
\fracd{\| \rho \|_2^2}{t} + 
\chi(A) e^{A/\tau}\int_\Omega \rho^2(\cdot,A,\cdot)\, dx
+2 \int_0^A&\int_\Omega (D+d) |\nabla\rho|^2 \, e^{a/\tau}dxda
\\
&+2 \int_0^A\int_\Omega \mu \rho^2 \, e^{a/\tau}dxda =
\int_\Omega\frac{\xi^2}{\tau^2} Q^2 \, dx
+ \frac{1}{\tau} \| \rho \|_2^2.
\end{aligned}\end{equation}
As the second, third and fourth terms of the left hand side of equation
(\ref{EsL2.9}) are non negative, we may deduce
\begin{equation}\label{EsL2.10}
\fracd{\| \rho \|_2^2}{t} \leq \frac{1}{\tau^2} \| Q \|_2^2
+  \frac{1}{\tau} \| \rho \|_2^2.
\end{equation}
In the same way, multiplying now (\ref{M1.a}) by $Q$ and integrating
yields
\begin{equation}\label{EsL2.11}
\fracd{\| Q \|_2^2}{t} =  2\int_\Omega \frac{1-\xi}{\tau}  Q^2 \,dx
+ 2\int_0^A\int_\Omega \mu\rho Q\, e^{a/\tau}dxda
 +2\chi(A) e^{A/\tau} \int_\Omega \rho(\cdot,A,\cdot)  Q\,dx.
\end{equation}
Concerning the second term of the right hand side of the last equality,
using Young's inequality and formula (\ref{EsL2.5}) and (\ref{EsL2.6}),
we get 
\begin{equation}\label{EsL2.12} \begin{aligned}
\int_0^A\int_\Omega \mu\rho Q\, e^{a/\tau}&dxda
\leq\int_0^A\ \bigg(\int_\Omega  \mu^2\rho^2 \, dx\bigg)^{1/2}
\bigg(\int_\Omega Q^2 \, dx\bigg)^{1/2} e^{a/\tau}da
\\
\leq
\tilde\mu \Bigg[
\int_0^t
\Big(\int_\Omega &\frac{\xi^2}{\tau^2} Q^2(t-a,\cdot)\, dx \Big)^{1/2}  
e^{a/\tau}da
\\
&
+H(A-t) e^{t/\tau} \int_t^A  
 \Big( \int_\Omega \rho_0^2(a-t,\cdot) e^{2(a-t)/\tau}\, dx \Big)^{1/2}  da
\Bigg] \|Q\|_2
\\
\leq \tilde\mu\max(e^{T/\tau},\frac{e^{T/\tau}}{\tau^2})
\Bigg[
\int_0^t
&\Big(\int_\Omega Q^2(t-a,\cdot)\, dx \Big)^{1/2}da
\\
&
+H(A-t) \int_t^A  
 \Big( \int_\Omega \rho_0^2(a-t,\cdot) e^{2(a-t)/\tau}\, dx \Big)^{1/2} da
\Bigg] \|Q\|_2,
\end{aligned}\end{equation}
where $H(a)=0$ if $a<0$ and $H(a)=1$ if $a\in[0,+\infty]$ and where
$\tilde \mu =\sup\{\mu(t,a,x),\, t\in [0,T), a\in [0,A), x\in\Omega \}.$
\newline
The third term of the right hand side of (\ref{EsL2.11}) may also be estimated:
\begin{equation}
\chi(A) e^{A/\tau}\int_\Omega \rho(\cdot,A,\cdot) Q \, dx \leq 
\chi(A) e^{A/\tau}\bigg( \int_\Omega \rho^2(\cdot,A,\cdot)\, dx\bigg)^{1/2}
\bigg(\int_\Omega Q^2 \, dx  \bigg)^{1/2}.
\end{equation}
Hence, applying again  (\ref{EsL2.5}) and  (\ref{EsL2.6}),
\begin{equation}\begin{aligned}\label{EsL2.13}
\chi(A) e^{A/\tau}\int_\Omega \rho(\cdot,A,\cdot) Q \, dx &\leq 
\chi(A) e^{A/\tau}\bigg(\int_\Omega Q^2(t-A,\cdot) \, dx  \bigg)^{1/2}
\|Q\|_2,\text{ ~ if } t>A,
\\
&
\leq \chi(A) e^{2t/\tau}
\bigg(\int_\Omega \rho_0^2(A-t,\cdot)  e^{2(A-t)/\tau}\, dx  \bigg)^{1/2}
\|Q\|_2, \text{ ~ if } t<A.
\end{aligned}\end{equation}
Using (\ref{EsL2.12}), (\ref{EsL2.13}) and (\ref{Hyprho0}) in 
(\ref{EsL2.11}), for four non negative constant $c_1$, $c_2$, $c_3$ and
 $c_4$ we have 
\begin{equation}\label{EsL2.14}
\fracd{\|Q\|_2^2}{t} \leq
\bigg( c_1\|Q\|_2  + 
\int_0^t \Big(\int_\Omega Q^2(t-a,\cdot) \, dx  \Big)^{1/2} da+
c_3\|\rho_0\|_2 +
c_4\Big(H(t-A)\int_\Omega Q^2(t-A,\cdot) \, dx  \Big)^{1/2} \bigg)\|Q\|_2.
\end{equation}
Setting 
\begin{equation}\label{EsL2.15}
F(t) = \sup_{s\in[0,T)} \|Q(s)\|_2^2 + \sup_{s\in[0,T)} \|\rho(s)\|_2^2,
\end{equation}
we have 
\begin{equation}\label{EsL2.16}
F'(t)\leq \max (0, \fracd{\|Q(t)\|_2^2}{t}) +
\max (0, \fracd{\|\rho(t)\|_2^2}{t}),
\end{equation}
and from (\ref{EsL2.10}) and (\ref{EsL2.14}) we deduce
\begin{equation}\label{EsL2.17}
F'(t)\leq \big((c_1+ c_2 T+ c_3+c_4)\sqrt{F(t)}\big)\sqrt{F(t)}+
\big(\frac{1}{\tau^2} +\frac{1}{\tau}\big)F(t)\leq c_5 F(t),
\end{equation}
for a constant $c_5$, from which we deduce that
\begin{equation}\label{EsL2.18}
F(t)\leq  e^{c_5T}F(0),
\end{equation}
and, as a consequence, that for a constant $c_6$,
\begin{equation}\label{EsL2.20}
\|Q(t)\|_2^2 + \|\rho(t)\|_2^2\leq c_6 (\|Q_0\|_2^2 + \|\rho_0\|_2^2). 
\end{equation}
Finally we get (\ref{EsL2}) as a consequence of (\ref{EsL2.20}), this
ends the proof. 
\end{proof}

As a consequence of the Lemma \ref{LemEs2} we have the following 
Corollary whose biological meaning is: Biomass cannot be created
from nothing.
\begin{cor}\label{LemEs3}
If the assumptions of Theorem \ref{Thm1} are realized and if 
$\rho_0=0$ and $Q_0=0$ then the solution $(Q,\rho)$ given by the 
Theorem satisfies 
\begin{equation}
\rho=0 \text { ~ ~ and ~ ~ } Q=0.
\end{equation}
\end{cor}
In order to establish the previous estimates, we have assumed that
$Q\geq 0$ and $\rho\geq 0$. We can show that this is a consequence of the 
non negativity of $Q_0$ and $\rho_0$.
\begin{lem}\label{LemEs5}
If the assumptions of Theorem \ref{Thm1} are true and if there exists
a solution $(Q,\rho)$ to system (\ref{M1.a})-(\ref{M1.f}), then
$Q\geq 0$ and $\rho\geq 0.$
\end{lem}
The proof of this Lemma is close to the one of Lemma \ref{LemEs2}.
\begin{proof}
We define
\begin{equation}\label{EsL2.23}
\rho^- =\min(\rho,0),  ~~~~~~~ Q^- = \min(Q,0).
\end{equation}
Of course,
\begin{equation}\label{EsL2.24}
\rho^-(0,\cdot,\cdot) =0,  ~~  Q^-(0,\cdot) = 0,  ~~
\rho^-(\cdot,0,\cdot) = \frac{\xi}{\tau} Q^-.
\end{equation}
Now, multiplying (\ref{EsL2.24}) by $e^{a/\tau}\rho^-$ and
integrating, since 
\begin{equation}\label{EsL2.25}
\fracp{\rho}{t}\rho^- = \fracp{\rho^-}{t}\rho^-,
\end{equation}
proceeding as while establishing (\ref{EsL2.10}), we get
\begin{equation}\label{EsL2.26}
\fracd{\| \rho^- \|_2^2}{t} \leq 
\frac{1}{\tau^2}\| Q^- \|_2^2 + \frac{1}{\tau}\| \rho^- \|_2^2.
\end{equation}
Multiplying now (\ref{M1.a}) by $Q^-$ and integrating gives
\begin{equation}\begin{aligned}
\fracd{\| Q^- \|_2^2}{t}  
&=2\int_\Omega \frac{1-\xi}{\tau} {Q^-}^2 \, dx 
+ 2\int_0^A\int_\Omega \mu\rho Q^- e^{a/\tau}  \, dxda 
+2 \chi(A) e^{A/\tau}\int_\Omega \rho(\cdot,A,\cdot) Q^- \,dx 
\\
&\leq 2\int_\Omega \frac{1-\xi}{\tau} {Q^-}^2 \, dx
+ 2\int_0^A\int_\Omega \mu\rho^- Q^- e^{a/\tau}  \, dxda 
+2 \chi(A) e^{A/\tau}\int_\Omega \rho^-(\cdot,A,\cdot) Q^- \,dx.
\end{aligned}\end{equation}
Since, as we  had (\ref{EsL2.5}) and (\ref{EsL2.6}), we have here:
\begin{equation}
\int_\Omega {\rho^-}^2 \,dx= \int_\Omega {\rho_0^-}^2 \,dx 
\text{ if } t<a \text { ~ and ~ } 
\int_\Omega {\rho^-}^2 \,dx= 
\int_\Omega \frac{\xi^2}{\tau^2}{Q^-}^2(t-a,\cdot) \,dx
\text{ if } t>a ,
\end{equation}
we can finish the proof as in the proof of Lemma \ref{LemEs2} and get
\begin{equation}
\| Q^-(t) \|+ \| \rho^-(t) \| \leq 
c'\big( \| Q^-_0 \|+  \| \rho^-_0 \|\big)=0,
\end{equation}
giving the Lemma.
\end{proof}
\subsection{Estimates on the derivatives}\label{SecDerEs}
As a by product of the proof of Lemma \ref{LemEs2}, we can deduce from
(\ref{EsL2.9}) the following result insuring a first control on the
regularity of $\rho$.
\begin{cor}\label{LemEs4}
If the assumptions and the conclusion of Theorem \ref{Thm1} are valid
then, for any $0\leq s\leq T$,
\begin{equation}\label{EsL2.22}
\int_0^s \| \nabla\rho\|_2^2 \,dt = 
\int_0^s\int_0^A\int_\Omega | \nabla\rho |_2^2 \,dxdadt
\leq \frac cd,
\end{equation}
for a constant $c$ (depending only on  $A$, $\tau$, $T$,  
$\tilde \mu$, $\| \rho_0 \|_{2}$ and $\| Q_0 \|_{2}$).
\end{cor}
\begin{proof}
Integrating (\ref{EsL2.9}) from 0 to $s$ yields:
\begin{multline}\label{EsL2.9.1}
\chi(A) e^{A/\tau}\int_0^s\int_\Omega \rho^2(\cdot,A,\cdot)\, dxdt+
2\int_0^s \int_0^A\int_\Omega (D+d) |\nabla\rho|^2 \, e^{a/\tau}dxdadt
\\
+2 \int_0^s\int_0^A\int_\Omega \mu \rho^2 \, e^{a/\tau}dxdadt =
\int_0^s\int_\Omega\frac{\xi^2}{\tau^2} Q^2 \, dxdt
+\frac{1}{\tau} \int_0^s \| \rho \|_2^2 \, dt
+\| \rho(0) \|_2^2 - \| \rho (s)\|_2^2.
\end{multline}
Using the previous estimate concerning $\rho$ and $Q$, we may deduce
\begin{equation}\label{EsL2.9.2}
d \int_0^s \int_0^A\int_\Omega |\nabla\rho|^2 \, e^{a/\tau}dxdadt\leq c,
\end{equation}
and (\ref{EsL2.22}) follows.
\end{proof}

~

Because of the form of the non linearity in (\ref{M1.a})-(\ref{M1.f}),
we need a supplementary estimate concerning 
\begin{equation}\label{Eslap1}
\int_0^T \| \Delta \rho \|_2^2 \,dt =
\int_0^T \int_0^A \int_\Omega| \Delta \rho |_2^2\,dxdadt.
\end{equation}
This estimate is a consequence of an estimate on $\| \nabla \rho \|_4$
and on
\begin{gather}\label{Eslap2}
\| \rho \|_\infty = 
\sup \big\{|\rho(t,a,x)|, t\in [0,T),a\in[0,A), x\in\Omega \big\}, ~~
\| Q \|_\infty = 
\sup \big\{|Q(t,x)|, t\in [0,T), x\in\Omega \big\}, ~~
\end{gather}
that we now set.
\begin{lem}\label{LemEs6}
The solution $(Q,\rho)$ given by Theorem \ref{Thm1} satisfies
\begin{gather}
\| \rho \|_\infty + \| Q \|_\infty \leq k,
\label{Eslap4}\end{gather}
where $k$ is a constant depending only on  
$T$, $A$, $\sup_{t\in[0,T)}\|\rho(t) \|_{2}$ and 
$\sup_{t\in[0,T)}\| Q(t) \|_{2}$
(which are estimated by Lemma \ref{EsL2})
\end{lem}
\begin{proof}
As we already see that $\|\rho(t) \|_{2}$ is bounded, using 
a method similar to Ladyzenskaja, Solonnikov and Ural'ceva \cite{LSU},
(paragraph III-8) we deduce that
\begin{gather}\label{Eslap5}
\big| \chi(A) \rho(\cdot,A,\cdot) e^{A/\tau} \big| \leq k_1,
\end{gather}
where the constant $k_1$ only depends only on $A$ and 
$\sup_{t\in[0,T)}\|\rho(t) \|_{2}$. Then defining
\begin{gather}\label{Eslap6}
\tilde P =\int_0^A\rho e^{a/\tau} \,da,
\end{gather}
we deduce from (\ref{M1.b}) that $\tilde P$ is solution to the following 
parabolic equation
\begin{gather}\label{Eslap6.1}
\fracp{\tilde P}{t} + \chi(A) \rho(\cdot,A,\cdot) e^{A/\tau} =
- \int_0^A \mu \rho e^{a/\tau} \,da
+\nabla\cdot \big((D+d)\nabla \tilde P \big) 
+\frac{\xi}{\tau}Q,
\end{gather}
from which we get that
\begin{gather}\label{Eslap7}
\|\tilde P \|_\infty = 
\sup\big\{ |\tilde P(t,x) |, t\in [0,T), x\in\Omega\big\}
\leq k_2,
\end{gather}
where $k_2$ only depends on $\ds\sup_{t\in[0,T)}\|\rho(t) \|_{2}$,
$\ds\sup_{t\in[0,T)}\| Q(t) \|_{2}$ and 
$\ds \sup_{x\in\Omega}|\tilde P_0| = 
\sup_{x\in\Omega}|\int_0^A \rho_0e^{a/\tau} \,da |$
which is finite by assumption.
\newline
Then, (\ref{Eslap7}) and (\ref{Eslap5}) give that $\| Q \|_{\infty}$
is finite, and as a consequence, (\ref{M1.b}) and (\ref{M1.c})
finally give the bound on $\| \rho \|_{\infty}$,
ending the proof.
\end{proof}
\begin{lem}\label{LemEs7}
The solution $(Q,\rho)$ given by Theorem \ref{Thm1} satisfies 
\begin{equation}\label{Eslap8}
\| \nabla Q \|_4 +
\| \nabla \rho \|_4 = 
\bigg( \int_\Omega | \nabla Q |^4\, dx\bigg)^{1/4}+
\bigg(\int_0^A \int_\Omega | \nabla \rho |^4 
\,e^{a/\tau}dxda \bigg)^{1/4}\leq \frac{C}{d},
\end{equation}
for a constant $C$ ($C$ does not depend on $d$).
\end{lem}
 \begin{proof}
Multiplying equation (\ref{M1.b}) by 
$-\nabla\cdot \big(|\nabla\rho|^2\nabla\rho  \big)e^{a/\tau}$, and 
integrating in $a$ and $x$ gives:
\begin{multline}\label{Eslap9}
-\int_0^A\int_\Omega\frac{\partial\rho}{\partial t}
\;\nabla\cdot \big(|\nabla\rho|^2\nabla\rho  \big)\,e^{a/\tau}dxda
-\int_0^A\int_\Omega\frac{\partial\rho}{\partial a}
\;\nabla\cdot \big(|\nabla\rho|^2\nabla\rho  \big)\,e^{a/\tau}dxda
=
\\
\int_0^A\int_\Omega\mu\rho
\;\nabla\cdot \big(|\nabla\rho|^2\nabla\rho  \big)\,e^{a/\tau}dxda
-\int_0^A\int_\Omega\nabla\cdot\big[\big(D+d\big)\nabla\rho\big]
\;\nabla\cdot \big(|\nabla\rho|^2\nabla\rho  \big)\,e^{a/\tau}dxda.
\end{multline}
Making  a double integration by part, and following a straightforward
computation procedure, we get
\begin{multline}\label{Eslap10}
\int_0^A\int_\Omega\nabla\cdot\big[\big(D+d\big)\nabla\rho\big]
\;\nabla\cdot \big(|\nabla\rho|^2\nabla\rho  \big)\,e^{a/\tau}dxda=
\\
\int_0^A\int_\Omega\big(D+d\big)
\big(
|\nabla\rho|^2|\nabla^2\rho|^2  + 2 \HH(\nabla\rho,\nabla^2\rho)
\big)\,e^{a/\tau}dxda +\EE,
\end{multline}
where 
\begin{gather}
\ds|\nabla^2\rho|^2=\big( \fracpsec{\rho}{x_1} \big)^2+
2\big(\fracpsecc{\rho}{x_1}{x_2}\big)^2+
\big(\fracpsec{\rho}{x_2}\big)^2,
\\
\ds\HH(\nabla\rho,\nabla^2\rho)=
\Big(\fracp{\rho}{x_1}\fracpsec{\rho}{x_1}+
\fracp{\rho}{x_2}\fracpsecc{\rho}{x_1}{x_2}
\Big)^2
+\Big(\fracp{\rho}{x_1}\fracpsecc{\rho}{x_1}{x_2}+
\fracp{\rho}{x_2}\fracpsec{\rho}{x_2}
\Big)^2,
\end{gather}
and
\begin{multline}\label{Eslap11}
\EE= 
\int_0^A \int_\Omega \sum_{i,j=1}^2 
\Big(\fracp{D({\cal M},Q,P)}{x_i} \fracp{\rho}{x_j}\Big)
\\
\Big(|\nabla\rho|^2 \fracpsecc{\rho}{x_i}{x_j} 
+2
\big( \fracp{\rho}{x_1}\fracpsecc{\rho}{x_1}{x_j} +
\fracp{\rho}{x_2}\fracpsecc{\rho}{x_j}{x_2}\big)\fracp{\rho}{x_i}\Big)
\,e^{a/\tau}dxda.
\end{multline}
Since 
\begin{gather}\label{Eslap11.1}
\fracp{D({\cal M},Q,P)}{x_i} =
\fracp{D}{{\cal M}}\fracp{{\cal M}}{x_i}
+\fracp{D}{Q} \fracp{Q}{x_i}
+\fracp{D}{P}\fracp{P}{x_i},
\end{gather}
in view of the regularity of $D$, of equation (\ref{Defm}) that gives a 
control on $\partial{\cal M} / \partial x_i$ in terms of
$\partial{P} / \partial x_i$, we get for a constant $C_1$
\begin{align}
|\EE| &\leq C_1\int_0^A \int_\Omega \big(|\nabla P| +|\nabla Q|\big)
|\nabla\rho|^3 \,|\nabla^2\rho| \,e^{a/\tau}dxda
\nonumber
\\
&\leq \frac{d}{4} \int_0^A \int_\Omega|\nabla\rho|^2 \,|\nabla^2\rho|^2
\,e^{a/\tau}dxda
+ \frac{C_1^2}{d} \int_0^A \int_\Omega \big(|\nabla P| +|\nabla Q|\big)^2
|\nabla\rho|^2\,e^{a/\tau}dxda;
\label{Eslap12}
\end{align}
in order to get the last expression in (\ref{Eslap12}), we used
$UV\leq\frac d4 U^2 + \frac 1d V^2$ with 
$V=C_1\big(|\nabla P| +|\nabla Q|\big)|\nabla\rho|$.
\newline
Concerning the other terms of (\ref{Eslap9}), since 
$\ds \fracp{|\nabla\rho|^4}{t}=2\fracp{|\nabla\rho|^2}{t}|\nabla\rho|^2=
4\fracp{\nabla\rho}{t}\cdot \nabla\rho \,|\nabla\rho|^2$,
making an integration by parts, we get
\begin{gather}\label{Eslap12bis}
-\int_0^A\int_\Omega\fracp{\rho}{t}
\;\nabla\cdot \big(|\nabla\rho|^2\nabla\rho  \big)\,e^{a/\tau}dxda=
\frac 14 \fracd{\|\nabla\rho\|_4^4}{t}.
\end{gather}
In a similar way,
\begin{multline}\label{Eslap13}
-\int_0^A\int_\Omega\frac{\partial\rho}{\partial a}
\;\nabla\cdot \big(|\nabla\rho|^2\nabla\rho  \big)\,e^{a/\tau}dxda=
\frac 14 \int_0^A\int_\Omega\fracp{(|\nabla\rho|^4\,e^{a/\tau})}{a}dxda
-\frac{1}{4\tau}\|\nabla\rho\|_4^4
\\
=\frac 14 \chi(A)\int_\Omega|\nabla\rho(\cdot,A,\cdot)|^4 \,e^{A/\tau}dx
-\frac 1{4\tau^4} \int_\Omega \big(\xi^4+(\fracp{\xi}{Q} Q)^4 \big)|\nabla Q|^4 \,dx
-\frac{1}{4\tau}\|\nabla\rho\|_4^4,
\end{multline}
and
\begin{multline}\label{Eslap14}
\bigg|
\int_0^A\int_\Omega\mu\rho
\;\nabla\cdot \big(|\nabla\rho|^2\nabla\rho  \big)\,e^{a/\tau}dxda
\bigg|
\\
=\bigg|
-\int_0^A\int_\Omega\mu |\nabla\rho|^4\,e^{a/\tau}dxda
-\int_0^A\int_\Omega\nabla\mu\, \rho\,|\nabla\rho|^2\,e^{a/\tau}dxda
\bigg|
\leq C_2 \big( \|\nabla\rho\|_4^4 +1  \big),
\end{multline}
using the regularity of $\mu$ and the estimate on $\sup(\rho)$ given by
Lemma \ref{LemEs6}.
\newline
The regularity of $\xi$ and the estimate on $\sup(Q)$ give
$|\xi^4+(\fracp{\xi}{Q} Q)^4 |\leq C_3$ for a constant $C_3$.
Hence (\ref{Eslap9}) yields
\begin{multline}\label{Eslap15}
\frac 14 \fracd{\|\nabla\rho\|_4^4}{t}
+\frac 14 \chi(A)\int_\Omega|\nabla\rho(\cdot,A,\cdot)|^4 \,e^{A/\tau}dx
+\int_0^A\int_\Omega\big(D+d\big)
\big(|\nabla\rho|^2 |\nabla^2\rho|^2 \big)\,e^{a/\tau}dxda
\\
+2 \int_0^A\int_\Omega\big(D+d\big)\,\HH(\nabla\rho,\nabla^2\rho)\,e^{a/\tau}dx
\leq
\frac{C_3}{4\tau^4}\int_\Omega|\nabla Q|^4 \,dx
+\frac{1}{4\tau}\|\nabla\rho\|_4^4
+C_2 \big( \|\nabla\rho\|_4^4 +1  \big)
\\
+\frac{d}{4} \int_0^A \int_\Omega|\nabla\rho|^2 \,|\nabla^2\rho|^2
\,e^{a/\tau}dxda
+\frac{C_1^2}{d} 
\bigg(\int_0^A \int_\Omega |\nabla\rho|^4\,e^{a/\tau}dxda
+C_4\big(\|\nabla Q\|_4^4 + \|\nabla\rho\|_4^4\big)
\bigg),
\end{multline}
and passing the fourth term of the right hand side in the left hand side
we can deduce
\begin{gather}\label{Eslap16}
\fracd{\|\nabla\rho\|_4^4}{t}\leq
\frac{C_5}{d} 
\big( \|\nabla\rho\|_4^4 + \|\nabla Q\|_4^4 +1\big),
\end{gather}
for a constant $C_5$.

Multiplying equation (\ref{M1.b}) by 
$-\nabla\cdot \big(|\nabla\rho|^2\nabla\rho \big)e^{4a/\tau}$, and 
making the same operations as previously, we obtain an inequality
which is (\ref{Eslap15}) with $e^{a/\tau}$ replaced by $e^{4a/\tau}$
and $\|\nabla Q\|_4^4 $ replaced by 
$\int_0^A \int_\Omega|\nabla\rho|^4\,e^{4a/\tau}dxda$. From
this, we can deduce
\begin{multline}\label{FM1}
\fracd{\big(\int_0^A \int_\Omega|\nabla\rho|^4\,e^{4a/\tau}dxda\big)}{t}
+ \chi(A) \int_\Omega |\nabla\rho(\cdot,A,\cdot)|^4\,e^{4A/\tau}dx
\\
\leq\frac{C_6}{d}\Big(
\int_0^A \int_\Omega|\nabla\rho|^4\,e^{4a/\tau}dxda+\|\nabla Q\|_4^4 +1
\Big).
\end{multline}

On another hand, computing the gradient of (\ref{M1.a}), and
multiplying by $\nabla Q |\nabla Q|^2$ yields
\begin{multline}\label{Eslap17}
\frac 14 \fracd{\|\nabla Q\|_4^4}{t}
=\frac 1\tau \int_\Omega\fracp{\xi}{Q} Q |\nabla Q|^4 \,dx 
+\int_\Omega\frac{1-\xi}{\tau}|\nabla Q|^4  \,dx 
+\int_\Omega \bigg(\int_0^A\nabla\mu \, \rho \,e^{a/\tau}da \bigg)
\cdot \nabla Q |\nabla Q|^2\,dx
\\
+\int_\Omega \bigg(\int_0^A\mu \nabla\rho \,e^{a/\tau}da \bigg)
\cdot \nabla Q |\nabla Q|^2\,dx
+ \chi(A) \int_\Omega\Big(\nabla\rho(\cdot,A,\cdot) e^{A/\tau}\Big)
\cdot \nabla Q |\nabla Q|^2\,dx.
\end{multline}
Because of the regularity of $\xi$ and $\mu$ and of Lemma \ref{LemEs6},
since, applying Young's inequality,
\begin{multline}\label{FM2}
\int_\Omega \bigg(\int_0^A\mu \nabla\rho \,e^{a/\tau}da \bigg)
\cdot \nabla Q |\nabla Q|^2\,dx 
\\
\leq
\Bigg(\int_\Omega \bigg(\int_0^A\mu |\nabla\rho| \,e^{a/\tau}da \bigg)^4
\,dx\Bigg)^{1/4}
\Bigg(\int_\Omega \bigg(|\nabla Q|^3\bigg)^{4/3}\,dx  
\Bigg)^{3/4}
\\
\leq \tilde\mu\Bigg(\int_0^A\int_\Omega|\nabla\rho|^4 
\,e^{4a/\tau}dxda \Bigg)^{1/4}
\|\nabla Q\|_4^3,
\end{multline}
and
\begin{multline}\label{FM3}
\chi(A)\int_\Omega\big( \nabla\rho(\cdot,A,\cdot) \,e^{A/\tau}\big)
\cdot \nabla Q |\nabla Q|^2\,dx
\leq\chi(A)\bigg(\int_\Omega|\nabla\rho(\cdot,A,\cdot)|^4 
\,e^{4A/\tau}\,dx\bigg)^{1/4}
\|\nabla Q\|_4^3
\\
\leq \chi(A)\frac{C'_6}{d} \bigg(\Big(
\int_0^A \int_\Omega|\nabla\rho|^4\,e^{4a/\tau}dxda\Big)^{1/4} + \|\nabla Q\|_4+1
\bigg)\|\nabla Q\|_4^3,
\end{multline}
we deduce from (\ref{Eslap17})
\begin{equation}\label{Eslap18}
\fracd{\|\nabla Q\|_4^4}{t} \leq C_4 (\|\nabla Q\|_4^4+1)
+\frac{C_5}{d}\bigg(\Big( \int_0^A \int_\Omega|\nabla\rho|^4\,e^{4a/\tau}dxda\Big)^{1/4} 
+ \|\nabla Q\|_4+1\bigg)\|\nabla Q\|_4^3.
\end{equation}
Inequalities (\ref{Eslap16}), (\ref{FM1}), (\ref{Eslap18}) and
the assumptions on $Q_0$ and $\rho_0$ give 
\begin{equation}\label{FM4}
\| \nabla Q \|_4 
+\| \nabla \rho \|_4 +
\Big( \int_0^A \int_\Omega|\nabla\rho|^4\,e^{4a/\tau}dxda\Big)^{1/4}
\leq \frac{C}{d},
\end{equation}
and finally the Lemma.
\end{proof}

\begin{lem}
The solution $(Q,\rho)$ given by Theorem \ref{Thm1} satisfies 
\begin{gather}
\| \nabla Q(t) \|^2_2 \leq \frac{c}{d}, ~~ 
\| \nabla \rho(t) \|^2_2 \leq \frac{c}{d}\text{ ~ for any } 0\leq t \leq T,
\label{Eslap19}
\\
\int_0^s\| \Delta\rho(t) \|^2_2\, dt  =\int_0^s\int_0^A\int_\Omega 
| \Delta\rho(t) |^2 \,e^{a/\tau} dxdadt \leq \frac{c}{d^2},
\text{ ~ for any } 0\leq s  \leq T,
\label{Eslap20}
\\
\int_0^s\bigg \| \fracp{Q}{t}(t) \bigg \|^2_2 dt \leq c,~~
 \int_0^s\bigg \| \fracp{\rho}{t}(t) \bigg \|^2_2 dt \leq \frac{c}{d^2}\text{ ~ and ~ }
\int_0^s\bigg \| \fracp{\rho}{a}(t) \bigg \|^2_2 dt\leq \frac{c}{d^2}
\text{ ~ for any } 0\leq s  \leq T,
\label{Eslap21}
\end{gather}
for a constant $c$ (which does not depend on $d$).
\end{lem}
\begin{proof}
Multiplying (\ref{M1.b}) by $(-\Delta\rho) e^{a/\tau}$ and integrating in $a$
and $x$ gives
\begin{multline}
\fracd{\| \nabla \rho\|_2^2}{t}
+ \chi(A) e^{A/\tau}\int_\Omega \big|\nabla\rho(\cdot,A,\cdot) \big|^2\,dx
+2\int_0^A\int_\Omega \mu \big|\nabla\rho \big|^2\,e^{a/\tau}dxda
\\
+2 \int_0^A\int_\Omega (D+d) \big|\Delta\rho \big|^2\,e^{a/\tau}dxda
= 
\int_0^A\int_\Omega \nabla\mu\cdot\nabla\rho \,\rho \,e^{a/\tau}dxda
\\
- \int_0^A\int_\Omega \nabla\big[D({\cal M},Q,P) \big] \cdot\nabla\rho 
\Delta\rho \,e^{a/\tau}dxda
+\frac 1\tau \|\nabla \rho \|_2^2 
+\int_\Omega \big(\xi^2 + (\fracp{\xi}{Q} Q)^2\big) |\nabla Q|^2\,dx
\\
\leq c_1 \big(\| \nabla\rho\|_2 +\| \nabla\rho\|^2_2 + \| \nabla Q\|^2_2\big) 
+ \int_0^A \int_\Omega \Big(\fracp{D}{{\cal M}} \nabla{\cal M}
+ \fracp{D}{P} \nabla P + \fracp{D}{Q} \nabla Q\Big) 
\nabla \rho \Delta\rho \,e^{a/\tau}dxda
\\
\leq c_1 \big( \| \nabla\rho\|_2 +\| \nabla\rho\|^2_2 + \| \nabla Q\|^2_2 \big)+
\frac 1d \int_0^A \int_\Omega 
\Big(\fracp{D}{{\cal M}} \nabla{\cal M}
+ \fracp{D}{P} \nabla P + \fracp{D}{Q} \nabla Q\Big)^2 
\big| \nabla \rho \big|^2\,e^{a/\tau}dxda
\\
+\frac d4  \int_0^A \int_\Omega \big| \Delta \rho \big|^2\,e^{a/\tau}dxda,
\end{multline}
for a constant $c_1$.
Now transferring the last term of the right hand side in the left
hand side and using the estimate on $\|\nabla\rho\|_4$ and $\|\nabla Q\|_4$  
which also give an estimate on $\|\nabla{\cal M}\|_4$ and $\|\nabla P\|_4$,
we get, for a constant $c_2$
\begin{gather}\label{CFT}
\fracd{\| \nabla \rho\|_2^2}{t} 
+\frac{5d}{4}  \| \Delta \rho\|_2^2 \leq c_2 
\big( \| \nabla \rho\|^2_2+\| \nabla Q\|^2_2 +\frac{1}{d}\big),
\end{gather}
In a similar way, we can also get 
\begin{gather}
\fracd{(\int_0^A\int_\Omega | \nabla\rho|^2 \,e^{2a/\tau}dxda)}{t}
\leq
c_3 \Big( \int_0^A\int_\Omega | \nabla\rho|^2 \,e^{2a/\tau}dxda + \|\nabla Q\|_2^2 
+\frac{1}{d}\Big),
\end{gather}
and 
\begin{gather}
\fracd{\| \nabla Q\|_2^2}{t}\leq
c_4\Big( \int_0^A\int_\Omega | \nabla\rho|^2 \,e^{2a/\tau}dxda + \|\nabla Q\|_2^2 \Big).
\end{gather}
From the three last inequalities we get (\ref{Eslap19}). 
Integrating (\ref{CFT}) from $0$ to $s$ gives  (\ref{Eslap20}). 
The estimate on $\| \nabla Q(t) \|_2$ is then obtained as the estimate 
on $\| \nabla Q(t) \|_4$.
Estimate (\ref{Eslap21})
is finally a direct consequence of equations (\ref{M1.a}) and (\ref{M1.b}).
\end{proof}
\begin{rem}
We could also prove that $\| \nabla Q \|_\infty$ and 
$\| \nabla \rho \|_\infty$ are bounded.
\end{rem}

\section{Existence and uniqueness of the solution}\label{SecExUn}
Once the a priori estimates are set, the proof of existence 
is classical and in the spirit of Ladyzenskaja, Solonnikov and Ural'ceva \cite{LSU}.
It essentially consists in linearization and passing to the limit using
the estimates.

In the following $L^p (\Omega)$ and $L^p ([0,A)\times\Omega,~e^{a/\tau}da dx)$
are the functional spaces associated with the norms (\ref{DefNorm}),
$W^{k,p}(\Omega)$ and $W^{k,p}([0,A)\times\Omega,~e^{a/\tau}da dx)$ are the
Sobolev spaces  composed of functions  whose derivatives up to order $k$
are in $L^p (\Omega)$ or  $L^p ([0,A)\times\Omega,~e^{a/\tau}da dx)$.
$L^p ([0,T)\times \Omega)$ and 
$L^p ([0,T)\times[0,A)\times\Omega,~e^{a/\tau}da dxdt)$
are the spaces of functions having finite norm
\begin{equation}\label{DefNorm33}
\bigg(\int_0^T\int_0^A \int_\Omega |\rho(t,a,x)|^p \,e^{a/\tau} dxdadt\bigg)^{1/p}
\text{ ~ or ~ }
\bigg(\int_0^T\int_\Omega |Q(t,x)|^p \, dxdt\bigg)^{1/p},
\end{equation}
and $W^{k,p}([0,T)\times\Omega)$ 
and $W^{k,p}([0,T)\times[0,A)\times\Omega,~e^{a/\tau}dadxdt)$
are their associated Sobolev spaces.
Finally, for a functional space $W$, $L^\infty(0,T;W)$ stands
for the functions whose norm in $W$ is finite for any $t\in[0,T]$.
\subsection{Linearization}\label{SecLin}
We linearize the system (\ref{M1.a})-(\ref{M1.f}) Then using classical results on
pde and ode, we give an existence and uniqueness result for the solution to
this linearized system.
\newline 
We set $Q^0=Q_0$ and $\rho^0=\rho_0$ and for $n\in\nit^*$, we consider
$(Q^n,\rho^n)$ solution to:
\begin{align}
& \frac{\partial Q^n}{\partial t}=\frac{1-\xi}{\tau}Q^n+
\int_0^{A}\rho^n(\cdot ,a,\cdot)e^{a/\tau}\mu(\cdot,a,\cdot)\,da +
\chi(A)\rho^n(\cdot,A,\cdot)e^{A/\tau},
~~\text{ on } [0,T) \times \Omega,
\label{LiM1.a}
\\
&\frac{\partial\rho^n}{\partial t}+\frac{\partial\rho^n}{\partial a}=
-\mu\rho^n+
\nabla\cdot\big[\big(D({\cal M}^{n-1},Q^{n-1},P^{n-1})+d\big)\nabla\rho^n\big],
~~\text{ on } [0,T) \times[0,A)\times \Omega,
\label{LiM1.b}
\\
&\rho^n(\cdot,0,\cdot)= \frac{\xi}{\tau}Q^{n},
~~\text{ on } [0,T) \times \Omega,
\label{LiM1.c}
\\
&\rho^n(0,\cdot,\cdot)= \rho_0,
~~\text{ on } [0,A) \times \Omega,
\label{LiM1.d}
\\
& Q^n(0,\cdot)=Q_0,
~~\text{ on } \Omega,
\label{LiM1.e}
\\
&\fracp{\rho^n}{\overrightarrow{\nu}} = 0,
~~\text{ on } [0,T) \times[0,A) \times \partial\Omega,
\label{LiM1.f}
\end{align}
where for $n\in\nit$,
\begin{equation}\label{LiDefP}
P^n(t,x) = \int_{a_{min}}^A \rho^n(t,a,x) e^{a/\tau} da,
\end{equation}
and ${\cal M}^n$ is solution to
\begin{align}
&\fracp{{\cal M}^n}{t} = \frac{1}{P_{max}-p_{max}} 
  ~H_r\bigg(\frac{P^n-p_{max}}{P_{max}-p_{max}}\bigg)
  ~H_r (1- {\cal M}^n) 
\nonumber
\\
&~~~~~~~~~~~~~~~~~~~~~~~~~~~~~~~~~~~~~~~~~~~~~~~~~~~~
- \frac{1}{p_{min}-P_{min}}
  ~H_r\bigg(\frac{p_{min}-P^n}{p_{min}-P_{min}}\bigg)~H_r ({\cal M}^n),
\nonumber
\\
& {\cal M}^n(0,\cdot)= {\cal M}_0.
\label{LiDefm}
\end{align}
\begin{thm}
Under assumptions (\ref{Hypmu}),(\ref{Hypxi}),(\ref{HypD}) and 
(\ref{Hypm0}), if 
$\rho_0\geq 0 \in L^1 \cap W^{2,2} \cap W^{1,4}
([0,A)\times\Omega,~e^{a/\tau}da dx)$ satisfies (\ref{Hyprho0})
and if $Q_0 \geq 0 \in L^1 \cap W^{2,2} \cap W^{1,4}(\Omega)$
then for any $n\in\nit$, there exists a unique solution
$(Q^n,\rho^n) \in C^0_b\big(0,T; (L^1 \cap W^{1,2}(\Omega) )
\times (L^1 \cap W^{1,2}([0,A)\times\Omega,~e^{a/\tau}da dx)) \big)$
to system (\ref{LiM1.a}) - (\ref{LiM1.f}) coupled with (\ref{LiDefP}) and 
(\ref{LiDefm}).
Moreover, $Q^n\geq 0$, $\rho^n\geq 0$ and $(Q^n,\rho^n)$ satisfies estimates
(\ref{EsL1}), (\ref{EsL2}), (\ref{Eslap4}), (\ref{Eslap8}),
(\ref{Eslap19}), (\ref{Eslap20}) and (\ref{Eslap21}) with constants
independent of $n$.
\end{thm}
The proof of this Theorem uses only classical pde and ode 
arguments. Hence we only sketch it.
\begin{proof}
The proof consists in an induction procedure. Because of the 
assumptions on $(Q_0,\rho_0)$  and the definition of $(Q^0,\rho^0)$,
the Theorem is true for $n=0$.
\newline
Then, if the Theorem is true for $n-1$, by regularization arguments,
we can get that that ${\cal M}^{n-1}$ exists and is unique on
$[0,T)\times\Omega$ and that 
${\cal M}^{n-1} \in C^0_b(0,T;W^{1,2}(\Omega)) \cap C^1_b(0,T;L^\infty(\Omega))$.
Hence we deduce that, for each $l\in\nit^*$, there exists a unique solution 
$(Q^{n,l},\rho^{n,l})\in C^0_b\big(0,T; (L^1 \cap W^{1,2}(\Omega) )
\times (L^1 \cap W^{1,2}([0,A)\times\Omega,~e^{a/\tau}da dx)) \big)$
to
\begin{align}
& \frac{\partial Q^{n,l}}{\partial t}=\frac{1}{\tau}Q^{n,l}-
\frac{\xi}{\tau}Q^{n,l-1}
+\int_0^{A}\rho^{n,l}(\cdot ,a,\cdot)e^{a/\tau}\mu(\cdot,a,\cdot)\,da +
\chi(A)\rho^{n,l}(\cdot,A,\cdot)e^{A/\tau},
\label{SwLiM1.a}
\\
&\frac{\partial\rho^{n,l}}{\partial t}+\frac{\partial\rho^n}{\partial a}=
-\mu\rho^{n,l}+
\nabla\cdot\big[\big(D({\cal M}^{n-1},Q^{n-1},P^{n-1})+d\big)\nabla\rho^{n,l}\big],
\label{SwLiM1.b}
\\
&\rho^{n,l}(\cdot,0,\cdot)= \frac{\xi}{\tau}Q^{n,l-1}, ~~ 
\rho^{n,l}(0,\cdot,\cdot)= \rho_0, ~~
\fracp{\rho^{n,l}}{\overrightarrow{\nu}}_{\big| \partial \Omega} = 0,
\label{SwLiM1.cdf}
\\
&
Q^{n,l}(0,\cdot)=Q_0,
\label{SwLiM1.e}
\end{align}
where $Q^{n,0}$ is defined as $Q^{n,0}=Q^{n-1}$.
This deduction involves first a classical semi-group or Galerkin
routine in order to deduce that there exists a unique solution 
$\rho^{n,l}$ to (\ref{SwLiM1.b}) - (\ref{SwLiM1.cdf}) as soon as 
the existence of $(Q^{n,l-1},\rho^{n,l-1})$ is achieved.
These routines are explained in Lions and Magenes \cite{LionsMagenes}, 
Ladyzenskaja, Solonnikov and Ural'ceva \cite{LSU}, or
--in a context close to our-- in Langlais \cite{Langlais1985}.
Once the existence of $\rho^{n,l}$ is established, 
the existence and uniqueness of $Q^{n,l}$ follows.
Now, following the way leading to (\ref{EsL2.10}) and (\ref{EsL2.14})
we deduce that $(Q^{n,l},\rho^{n,l})$ satisfies
\begin{gather}\label{EsL2.10.Sw}
\fracd{\| \rho^{n,l} \|_2^2}{t} \leq \frac{1}{\tau^2} \| Q^{n,l-1} \|_2^2
+  \frac{1}{\tau} \| \rho^{n,l} \|_2^2,
\\
\fracd{\|Q^{n,l}\|_2^2}{t} \leq
\bigg( c_1\|Q^{n,l}\|_2  + c'_1\|Q^{n,l-1}\|_2  
\int_0^t \Big(\int_\Omega (Q^{n,l})^2(t-a,\cdot) \, dx  \Big)^{1/2} da+
c_3\|\rho_0\|_2 
\nonumber
\\ 
~~~~~~~~~~~~~~~~~~~~~~~~~~~~~~~~~~~~~~~~~~~~~
+c_4\Big(H(t-A)\int_\Omega (Q^{n,l})^2(t-A,\cdot) \, dx  \Big)^{1/2} \bigg)\|
Q^{n,l}\|_2,
\label{EsL2.14.Sw}
\end{gather}
which is enough to deduce that $(\| \rho^{n,l} \|_2+\|Q^{n,l}\|_2) $
is bounded. As a consequence of this bound, we get that, for a subsequence
still denoted $l$,
$(Q^{n,l},\rho^{n,l}) \longrightarrow (Q^{n},\rho^{n})$ in 
$L^\infty\big(0,T; (L^2 (\Omega) )
\times (L^2 ([0,A)\times\Omega,~e^{a/\tau}da dx))$ weakly$-*$,
where $(Q^{n},\rho^{n})$ is solution to (\ref{LiM1.a}) - (\ref{LiM1.f}).
Finally, the estimates are led in the same way as in section 
\ref{SecAPrEst}. 
The uniqueness follows directly from the linear character of 
(\ref{LiM1.a}) - (\ref{LiM1.f}).
Hence, the Theorem is true for $n$.
\newline
The induction procedure is then straightforward to end the proof of 
the Theorem.
\end{proof}
\subsection{Existence}\label{SecEx}
From  estimates (\ref{Eslap19}), (\ref{Eslap20}) and (\ref{Eslap21})
we deduce that the sequence $(Q^n,\rho^n)$ is bounded in
$W^{1,2}([0,T)\times\Omega)\times W^{1,2}([0,T)\times[0,A)\times\Omega)$.
Hence, up to a subsequence still denoted $n$, we have
 $(Q^n,\rho^n)\longrightarrow (Q,\rho)$ in 
$W^{1,2}([0,T)\times\Omega)\times W^{1,2}([0,T)\times[0,A)\times\Omega,~e^{a/\tau}da dxdt)$
weakly, and then, in 
$L^{2}([0,T)\times\Omega)\times L^{2}([0,T)\times[0,A)\times\Omega,~e^{a/\tau}da dxdt)$
strongly.
\newline
From this we can also deduce that $P^n\longrightarrow P$ in
$L^{2}([0,T)\times\Omega)$ strongly, with $P$ defined from $\rho$ by 
(\ref{DefP}).
In view of (\ref{LiDefm}), we can deduce that $({\cal M}^n)$,
$(\partial{\cal M}^n/\partial t)$ and, taking the gradient of 
(\ref{LiDefm}), $(\nabla{\cal M}^n)$  are bounded in 
$L^{2}([0,T)\times\Omega)$. Then extracting again a subsequence
still denoted $n$, we deduce ${\cal M}^n\longrightarrow {\cal M}$
strongly, where ${\cal M}$ is solution to (\ref{Defm}).
\newline
Using now the regularity of $D$, we obtain 
$D({\cal M}^{n-1},Q^{n-1},P^{n-1})\longrightarrow
D({\cal M},Q,P)$ in $L^{2}([0,T)\times\Omega)$ strongly.
\newline
The regularity of trace operators gives
$\chi(A)\rho^{n}(\cdot,A,\cdot)\longrightarrow \chi(A)\rho(\cdot,A,\cdot)$,
$\rho^n(\cdot,0,\cdot)\longrightarrow \rho(\cdot,0,\cdot)$,
$\rho^n(0,\cdot,\cdot)\longrightarrow \rho(0,\cdot,\cdot)$,
$Q^n(0,\cdot)\longrightarrow Q^n(0,\cdot)$ weakly,
and using (\ref{Eslap20}),
$\partial \rho^n/\partial\overrightarrow{\nu}_{|\partial\Omega}
\longrightarrow 
\partial \rho/\partial\overrightarrow{\nu}_{|\partial\Omega}$ 
weakly.
\newline
Then passing to the limit in (\ref{LiM1.a}) - (\ref{LiDefm}),
we obtain that $(Q,\rho)$ is solution to 
(\ref{M1.a}) - (\ref{M1.f}) coupled with (\ref{DefP}) and 
(\ref{Defm}).

Once this existence result is established, 
using regularizations and truncations, we can start the
computations of section \ref{SecAPrEst} giving the 
additional regularity and the non negativity of
$Q$ and $\rho$.

It now remains to prove the uniqueness of the solution.

\subsection{Uniqueness}\label{SecUn}
Consider $(Q,\rho)$ with associated $P$ and ${\cal M}$ and 
$(\hat Q,\hat\rho)$ with associated $\hat P$ and $\hat {\cal M}$
two solutions of (\ref{M1.a}) - (\ref{M1.f}).
They both satisfy the estimates and the difference 
$(\tilde Q,\tilde\rho)=(Q-\hat Q,\rho-\hat\rho)$ satisfies
\begin{align}
& \frac{\partial \tilde Q}{\partial t}=\frac{1-\xi}{\tau}\tilde Q+
\int_0^{A}\tilde \rho(\cdot ,a,\cdot)e^{a/\tau}\mu(\cdot,a,\cdot)\,da +
\chi(A)\tilde \rho(\cdot,A,\cdot)e^{A/\tau},
\label{UqM1M1.a}
\\
&\frac{\partial\tilde \rho}{\partial t}+\frac{\partial\tilde \rho}{\partial a}=
-\mu\tilde \rho+
\nabla\cdot\big[\big(D({\cal M},Q,P)+d\big)\nabla\tilde \rho\big]-
\nabla\cdot\big[\big(D(\hat {\cal M},\hat Q,\hat P)-D({\cal M},Q,P)\big)
\nabla\hat\rho\big],
\label{UqM1M1.b}
\\ 
&\tilde\rho(\cdot,0,\cdot)= \frac{\xi}{\tau}\tilde Q, ~~~~~~~~~~~
\tilde\rho(0,\cdot,\cdot)= 0, ~~~~~~~~~~~
\tilde Q(0,\cdot)=0, ~~~~~~~~~~~
{\fracp{\tilde\rho}{\overrightarrow{\nu}}}_{|\partial\Omega} = 0.
\label{UqM1.cdef}
\end{align}
Multiplying (\ref{UqM1M1.b}) by $\tilde \rho e^{a/\tau}$ and integrating
gives
\begin{multline}\label{UqEsL2.9}
\fracd{\| \tilde \rho \|_2^2}{t} + 
\chi(A) e^{A/\tau}\int_\Omega \tilde \rho^2(\cdot,A,\cdot)\, dx
+2 \int_0^A\int_\Omega (D(\tilde{\cal M},\tilde Q,\tilde P)+d) 
|\nabla\tilde \rho|^2 \, e^{a/\tau}dxda
\\
= -2 \int_0^A\int_\Omega \mu \tilde \rho^2 \, e^{a/\tau}dxda 
+\int_\Omega\frac{\xi^2}{\tau^2} \tilde Q^2 \, dx
+ \frac{1}{\tau} \| \tilde \rho \|_2^2~~~~~~~~~~~~~~
\\~~~~~~~~~~~~~~~~~~~~~~~~
+ \int_0^A\int_\Omega \big( D(\hat{\cal M},\hat Q,\hat P)-D({\cal M},Q,P) \big)
\nabla\hat\rho \cdot \nabla\tilde\rho\, e^{a/\tau}dxda
\\
\leq k_1 \big( \| \tilde \rho \|_2^2 +\| \tilde Q \|_2^2\big) + 
k_2 \|\nabla\hat\rho\|_\infty\big( \| \tilde \rho \|_2 
+\| \tilde Q \|_2\big) \| \nabla\tilde \rho \|_2
\\
\leq k_1 \big( \| \tilde \rho \|_2^2 +\| \tilde Q \|_2^2\big) + 
\frac{k_2^2}{d}\|\nabla\hat\rho\|^2_\infty\big( \| \tilde \rho \|_2 +\| \tilde Q \|_2\big)^2+
\frac{d}{4}\| \nabla\tilde \rho \|_2^2,
\end{multline}
for constants $k_1$ and $k_2$.  Passing the term 
$\frac{d}{4}\| \nabla\tilde \rho \|_2^2$ in the left hand side gives,
for a constant $k_3$ 
\begin{equation}
\fracd{\| \tilde \rho \|_2^2}{t} \leq  \frac{k_3}{d} 
\big( \| \tilde \rho \|_2^2 +\| \tilde Q \|_2^2\big).
\end{equation}
Making the same, but multiplying (\ref{UqM1M1.b}) by 
$\tilde \rho e^{2a/\tau}$ yields
\begin{equation}
\fracd{\big(\int_0^A\int_\Omega | \tilde \rho |^2 \, e^{2a/\tau}dxda\big)}{t} 
+\chi(A) e^{2A/\tau} \int_\Omega  \tilde\rho(\cdot,A,\cdot) \,dx
\leq  \frac{k_3}{d} 
\bigg(\int_0^A\int_\Omega | \tilde \rho |^2 \, e^{2a/\tau}dxda +
\| \tilde Q \|_2^2\bigg).
\end{equation}
Multiplying (\ref{UqM1M1.a}) by $\tilde Q$ gives
\begin{equation}
\fracd{\| \tilde Q \|_2^2}{t} \leq\frac{k_4}{d} \bigg(
\int_0^A\int_\Omega | \tilde \rho |^2 \, e^{2a/\tau}dxda + \| \tilde Q \|_2^2
\bigg),
\end{equation}
for a constant $k_4$
From the last three inequalities we deduce, for a constant $K$
\begin{equation} 
\| \tilde Q \|_2^2+ \| \tilde \rho \|_2^2 + 
\int_0^A\int_\Omega | \tilde \rho |^2 \, e^{2a/\tau}dxda 
\leq 
K  \bigg( \| \tilde Q_{|t=0} \|_2^2+ \| \tilde \rho_{|t=0} \|_2^2 + 
\int_0^A\int_\Omega | \tilde \rho_{|t=0} |^2 \, e^{2a/\tau}dxda \bigg)
=0,
\end{equation}
giving $\tilde Q=\tilde \rho=0$ and then the uniqueness of the solution
to (\ref{M1.a}) - (\ref{M1.f}).
~~~~\hfill\rule{2mm}{2mm}\vskip3mm 
\bibliographystyle{plain}
\bibliography{biblio}

\begin{thebibliography}{10}

\bibitem{Andr1989}
V.~Andreasen.
\newblock Diesease regulation of age-structured host populations.
\newblock {\em Theo. Pop. Biol.}, 36(2):214--239, 1989.

\bibitem{Andr1992}
V.~Andreasen.
\newblock The effect of age-dependent host mortality on the dynamics of an
  endemic disease.
\newblock {\em Math. BioSciences}, 114:29--58, 1993.

\bibitem{Andr1995}
V.~Andreasen.
\newblock Instability in {SIR}-model with age-dependent susceptibility.
\newblock {\em Math. Pop. Dyn.}, 1:3--14, 1995.

\bibitem{AyDu}
B.~P. Ayati and T.~F. Dupont.
\newblock Galerkin methods in age and space for a population model with
  nonlinear diffusion.
\newblock {\em SIAM J. Numer. Anal.}, 40(3):1064--1076, 2002.

\bibitem{BusIann1983}
S.~Busenberg and M.~Iannelli.
\newblock A class of nonlinear diffusion problems in age-dependent population
  dynamics.
\newblock {\em Nonlinear Analysis, Theo., Meth., \& Appl.}, 7(5):501--529,
  1983.

\bibitem{DiBla1979}
G.~Di~Blasio.
\newblock Non linear age-dependent population diffusion.
\newblock {\em J. Math. Biol.}, 8:265--284, 1979.

\bibitem{DiBlaLam1978}
G.~Di~Blasio and L.~Lamberti.
\newblock An initial-boundary value problem for age-dependent population
  diffusion.
\newblock {\em SIAM J. Appl. Math.}, 35(3):593--616, November 1978.

\bibitem{EsiSha}
S.~E. Esipov and J.~A. Shapiro.
\newblock Kinetic model of \textit{{P}roteus mirabilis} swarm colony
  development.
\newblock {\em J. Math. Biol.}, 36:249--268, 1998.

\bibitem{Gueetal2001}
M.~Gu\'e, V.~Dupont, A.~Dufour, and O.~Sire.
\newblock Bacterial swarming: A biological time-resolved {FTIR-ATR} study of
  \textit{{P}roteus mirabilis} swarm-cell differentiation.
\newblock {\em BioChemistry}, 40:11938--11945, 2001.

\bibitem{Gur}
M.~E. Gurtin.
\newblock A system of equations for age-dependent population diffusion.
\newblock {\em J. theor. Biol}, 40:389--392, 1973.

\bibitem{GurMcCa1974}
M.~E. Gurtin and R.~C. Mac~Camy.
\newblock Non-linear age-dependent population dynamics.
\newblock 54:281--300, 1974.

\bibitem{GurMcCa}
M.~E. Gurtin and R.~C. Mac~Camy.
\newblock Diffusion model for age-structured population.
\newblock {\em Math. BioSciences}, 54:49--59, 1981.

\bibitem{Hua1994}
C.~Huang.
\newblock An age-dependent population model with nonlinear diffusion in
  {R}$^n$.
\newblock {\em Quat. of Appl. Math.}, LII(2):377--398, june 1994.

\bibitem{Kim1996}
M-Y. Kim.
\newblock Galerkin methods for a model of population dynamics with nonlinear
  diffusion.
\newblock {\em Num. Meth. for PDE}, 12:59--73, 1996.

\bibitem{KuLang1991}
M.~Kubo and M.~Langlais.
\newblock Periodic solutions for population dynamics problem with age-dependent
  and spatial structure.
\newblock {\em J. Math. Biol.}, 29:363--378, 1991.

\bibitem{LSU}
O.~A. Ladyzenskaja, Solonnikov~V. A., and Ural'ceva~N. N.
\newblock {\em Linear and Quasi-linear Equation of Parabolic Type}.
\newblock AMS, Translation of Mathematical Monographs, vol. 23.

\bibitem{Langlais1985}
M.~Langlais.
\newblock A nonlinear problem in age-dependent population diffusion.
\newblock {\em Siam J. Math. Anal.}, 16(3):510--529, 1985.

\bibitem{Langlais1988}
M.~Langlais.
\newblock Large time behavior in a nonlinear age-dependent population dynamics
  problem with spatial diffusion.
\newblock {\em J. Math. Biol.}, 26:319--346, 1988.

\bibitem{LionsMagenes}
J.~L. Lions and Magenes E.
\newblock {\em Probl\`emes aux limites non homog\`enes et applications}, volume
  17 , 18 of {\em Travaux et recherches math\'ematiques}.
\newblock Dunod.

\bibitem{LopTri1985}
L.~Lopez and D.~Trigiante.
\newblock A finite difference scheme for stiff problem arising in the numerical
  solution of a population dynamic model with spatial diffusion.
\newblock {\em Nonlinear Analysis, Theo., Meth., \& Appl.}, 9(1):1--12, 1985.

\bibitem{McCa}
R.~C. Mac~Camy.
\newblock A population model with nonlinear diffusion.
\newblock {\em J. Diff. Equ.}, 39:52--72, 1981.

\bibitem{Mar1981}
P.~Marcati.
\newblock Asymptotic behavior in age-dependent population dynamics with
  heredity renewal law.
\newblock {\em SIAM J. Math. Anal.}, 12(6):904--916, November 1981.

\bibitem{MedKaKo}
G.~E. Medvedev, T.~J. Kapper, and Koppel N.
\newblock A reaction-diffusion system with periodic front dynamics.
\newblock {\em SIAM J. Appl. Math.}, 60(5):1601--1638, 2000.

\bibitem{Mil1990}
F.~A. Milner.
\newblock A numerical method for a model of population dynamics with spatial
  diffusion.
\newblock {\em Comp. Math. Applic.}, 19(4):31--43, 1990.

\bibitem{Rauetal1996}
O.~Rauprich, M.~Matsuchita, C.~J. Weijer, F.~Siegert, S.~E. Esipov, and J.~A.
  Shapiro.
\newblock Periodic phenomena in \textit{{P}roteus mirabilis} swarm colony
  development.
\newblock {\em Jour. of Bacteriology.}, pages 6525--6538, Nov. 1996.

\end{thebibliography}
\end{document}